\documentclass[revision]{FPSAC2019}
\usepackage{lipsum}
\usepackage{amsmath}
\usepackage{bbm}
\usepackage{tikz}
\usepackage{pgf}
\usepackage{graphicx}
\usepackage{float}

\usepackage{todonotes}

\DeclareMathOperator{\st}{\,|\,}

\DeclareMathOperator{\C}{\mathbb{C}}
\DeclareMathOperator{\R}{\mathbb{R}}
\DeclareMathOperator{\Z}{\mathbb{Z}}

\DeclareMathOperator{\cyc}{cyc}

\def\s{{\mathfrak S}}
\def\D{{\mathcal D}}

\newtheorem{theorem}{Theorem}
\numberwithin{theorem}{section}
\newtheorem{lemma}[theorem]{Lemma}

\theoremstyle{definition}

\newtheorem{example}[theorem]{Example}

\title[Homogenized Linial Arrangement]{On the Homogenized Linial Arrangement: Intersection Lattice and Genocchi Numbers}

\author[Lazar \& Wachs]{Alexander Lazar\thanks{\href{mailto:alazar@math.miami.edu}{alazar@math.miami.edu}.} \and Michelle L. Wachs\thanks{Supported in part by NSF grant DMS 1502606; \href{mailto:wachs@math.miami.edu}{wachs@math.miami.edu}.}}

\address{Department of Mathematics, University of Miami, Coral Gables, FL 33124}

\received{November, 2018}

\revised{April 1, 2019}

\abstract{Hetyei recently introduced a hyperplane arrangement (called the homogenized Linial arrangement) and used the finite field method of Athanasiadis to show that its number of regions is a median Genocchi number. These numbers count a class of permutations known as Dumont derangements. Here, we take a different approach, which makes direct use of Zaslavsky's formula relating the intersection lattice of this arrangement to the number of regions. We refine Hetyei's result by obtaining a combinatorial interpretation of the Möbius function of this lattice in terms of variants of the Dumont permutations.  This enables us to derive a formula for the generating function of  the characterisitic polynomial of the arrangement.  The M\"obius invariant of the lattice turns out to be a (nonmedian) Genocchi number.  Our techniques also yield type B, and more generally Dowling arrangement, analogs of these results.}


\begin{document}

\maketitle
%% note that you DO NOT have to put your abstract here -- it is generated by \maketitle and the \abstract and \resume commands above
\section{Introduction}
Let $n \ge 1$. The {\em braid arrangement} is the hyperplane arrangement in $\R^n$ defined by
$$\mathcal A_{n-1} := \{ x_i - x_j = 0 : 1 \le i < j \le n \}.$$  Note that the hyperplanes of $\mathcal A_{n-1}$ divide $\R^n$ into $n!$ open cones of the form $$R_\sigma := \{{\bf x} \in \R^n : x_{\sigma(1)} < x_{\sigma(2)}< \dots < x_{\sigma(n)}\},$$
where $\sigma$ is a permutation in the symmetric group $\s_n$.

  A classical formula of Zaslavsky \cite{Facing_Up_Arrangements} gives the number of regions of any real hyperplane arrangement  $\mathcal A$ in terms of the M\"obius function of the intersection (semi)lattice $\mathcal L({\mathcal A})$ (which consists of intersections of collections of hyperplanes in $\mathcal{A}$, viewed as affine subspaces of $\R^n$, ordered by reverse containment).  Indeed, given any finite, ranked poset $P$ of length $r$, with a minimum element $\hat 0$, define the {\it characteristic polynomial} of $P$ to be 
 \begin{equation} \chi_P(t)  :=  \sum_{x \in P}\mu_P(\hat{0},x)t^{r-\text{rk}(x)},\end{equation}
  where $\mu_P(x,y)$ is the M\"obius function of $P$ and $\text{rk}(x)$ is the rank of $x$.  Zaslavsky's formula is
  \begin{equation} \label{zaseq} \#\{\text{regions of }\mathcal{A}\} = |\chi_{\mathcal L(\mathcal{A})}(-1)|. \end{equation}

  It is well known and easy to see that the lattice of intersections of the braid arrangement $\mathcal A_{n-1}$ is isomorphic to  the  lattice $\Pi_n$ of partitions of the set $[n]:=\{1,2\dots,n\}$.  It is also well known that the characteristic polynomial of $\Pi_n$ is given by
\begin{equation} \label{stireq} \chi_{\Pi_n} (t) = \sum_{k=1}^{n} s(n,k) t^{k-1},\end{equation}
  where $s(n,k)$ is the Stirling number of the first kind, which is equal to
  $(-1)^{n-k}$  times the number of permutations in $\s_n$ with exactly $k$ cycles; see \cite[Example 3.10.4]{EC1}. Hence
  $\chi_{\Pi_n}(-1) = (-1)^{n-1} |\s_n|$.  Therefore  from (\ref{zaseq}), we recover the result observed above that the number of regions of $\mathcal A_{n-1}$ is $n!$.
  
  In this extended abstract of \cite{Full_Paper}, we obtain analogous results for a  hyperplane arrangement introduced by Hetyei \cite{Alternation_Acyclic}. 
The {\em homogenized Linial arrangement}  is the hyperplane arrangement in $\R^{2n}$ defined, for $n \ge 2$, by\footnote{To justify our our indexing, we note that the length of the intersection lattice is $2n-3$.}
$$\mathcal{H}_{2n-3}:= \{x_i - x_j = y_i \st 1 \leq i < j \leq n\}.$$

Note that by intersecting $\mathcal{H}_{2n-3}$ with the subspace $y_1=y_2 =\dots = y_n =0$ one gets the braid arrangement $\mathcal A_{n-1}$.  Similarly by intersecting  $\mathcal{H}_{2n-3}$ with the subspace $y_1=y_2 =\dots = y_n =1$, one gets the Linial arrangement in $\R^n$,
$$ \{x_i - x_j = 1 \st 1 \leq i < j \leq n\}.$$ Postnikov and Stanley \cite{Defcox} show that the number of regions of the Linial arrangement is equal to the number of alternating trees on node set $[n+1]$. 

 Using the finite field method of Athanasiadis \cite{Finite_Field_Method}, Hetyei \cite{Alternation_Acyclic} obtains a recurrence for $\chi_{\mathcal L(\mathcal{H}_{2n-1})}(t)$ and uses it to show that 
\begin{equation} \label{heteq} |\chi_{\mathcal L(\mathcal{H}_{2n-1})}(-1) |=   h_{n},\end{equation}
where $h_n$ is a median Genocchi number.\footnote{In the literature the median Genocchi number $h_n$ is usually denoted $H_{2n+1}$.} 
Barsky and Dumont \cite[Theorem 1]{Barsky_Dumont} obtain the following generating function for the median Genocchi numbers
\begin{equation} \label{BDeq} \sum_{n\geq 1}h_{n}x^n = \sum_{n\geq 1}\frac{n!(n+1)!x^{n}}{\prod_{k=1}^n(1+k(k+1)x)}.\end{equation}

The median Genocchi numbers  also have numerous combinatorial interpretations.    One  of these interpretations is given in terms of a class of permutations called Dumont permutations; see \cite{Interpretations_Combinatoires} and \cite[Corollary 2.4]{Derangements_Genocchi}. Another is given in terms of surjective pistols in \cite[Corollary 2.2]{Derangements_Genocchi}. 

Here, we  study the intersection lattice $\mathcal L(\mathcal{H}_{2n-1})$.  We refine Hetyei's result (\ref{heteq}) by deriving a combinatorial formula for the M\"obius function of $\mathcal L(\mathcal{H}_{2n-1})$ in terms of permutations in $\s_{2n}$  similar to Dumont permutations, which we call D-permutations.  A key step in our proof is to show that $\mathcal L(\mathcal{H}_{2n-1})$  is isomorphic to the bond lattice of a certain bipartite graph. This bond lattice has  a nice description as  the induced subposet   of the partition lattice $\Pi_{2n}$ consisting of partitions all of whose nonsingleton blocks have odd smallest element and even largest element.     

Our M\"obius function result yields a combinatorial formula for the characteristic polynomial of $\mathcal L(\mathcal{H}_{2n-1})$ analogous to (\ref{stireq}) with  $\s_n$ replaced by the D-permutations   in $\s_{2n}$.
 By constructing a bijection between the D-permutations and surjective pistols, we recover Hetyei's result that $|\chi_{\mathcal L (\mathcal H_{2n-1})}(-1)|$ is a median Genocchi number. Moreover, we obtain the new result that the (nonmedian) Genocchi number\footnote{These are the signless Genocchi numbers; $g_n$ is usually denoted $(-1)^{n+1}G_{2n}$ in the literature.} $g_{n}$  is equal to  $|\mu_{\mathcal L(\mathcal{H}_{2n-1})}(\hat 0, \hat 1)| $, where
 $\hat 0$ and $\hat 1$ are the minimum and maximum elements of $\mathcal L (\mathcal H_{2n-1})$, respectively.  The bijection also enables us to derive a formula for the generating function $\sum_{n\ge 1} \chi_{\mathcal L (\mathcal H_{2n-1})}(t) x^n$, which reduces to the Barsky-Dumont formula (\ref{BDeq}) when $t= -1$ and to a similar formula of Barsky and Dumont in \cite{Barsky_Dumont} for  $\sum_{n\ge 1} g_n x^n$ when $t=0$.
 
 Our techniques also yield a type B analog of Hetyei's result and more generally a Dowling arrangement analog. 
 We define the {\em type B homogenized Linial arrangement} to be the hyperplane arrangement in $\R^{2n}$ defined by
  \begin{equation}  \mathcal{H}^B_{2n-1} =\{ x_i \pm x_j = y_i : 1 \le i < j \le n \} \cup \{x_i =y_i: i = 1\dots, n\}.\end{equation}
  We show that that the intersection lattice  of $\mathcal{H}^B_{2n-1}$ is isomorphic to an  induced subposet of the signed partition lattice $\Pi_{2n-1}^B$ and obtain  results for the M\"obius function and characteristic polynomial analogous to those for $\mathcal L(\mathcal{H}_{2n-1})$.  We use these results to prove the following generating function formula for the number of regions $r_n^B$ of $\mathcal{H}^B_{2n-1}$,
   \begin{equation} \label{BBDeq} \sum_{n\geq 1}r^B_n x^{n} = \sum_{n\geq 1}\frac{(2n)!x^n}{\prod_{k=1}^n( 1 + 2k(2k+1)x)},\end{equation}
 thereby providing a type B analog of (\ref{BDeq}).  We also obtain a type B analog of  the Barsky-Dumont formula  for the  generating function of the Genocchi numbers.

 Let $\omega$ be the primitive $m$th root of unity $e^{\frac{2\pi i }{m}}$.  For $m,n \ge 1$, the {\it Dowling arrangement}  is a hyperplane arrangement in $\C^n$ defined by
 \begin{equation} \label{doweq}  \{x_i -\omega^l x_j = 0: 1 \le i < j \le n, \, 0 \le l < m\} \cup \{x_i =0 : 1 \le i \le n \}.\end{equation}
 This is called a Dowling arrangement because its intersection lattice
 is isomorphic to the classical Dowling lattice $Q_n(\Z_m)$, which  is isomorphic to $\Pi_{n+1}$ when $m=1$, and to  $\Pi^B_n$ when $m=2$.  By introducing a Dowling analog of the homogenized Linial arrangement, we obtain unifying generalizations of the  types A and B results discussed above.  These generalizations include a polynomial analog of the formula $g_{n}=|\mu_{\mathcal L(\mathcal{H}_{2n-1})}(\hat 0, \hat 1)| $ involving a  polynomial analog of the Genocchi numbers known as the Gandhi polynomials.

\section{Preliminaries}

{\bf  Hyperplane Arrangements.} 
Let $k$ be a field (here $k$ is $\R$ or $\C$). A {\em hyperplane arrangement} $\mathcal{A} \subseteq k^n$ is a finite collection of affine codimension-$1$ subspaces of $k^n$.
The {\em intersection poset} of $\mathcal{A}$ is the poset $\mathcal{L}(\mathcal{A})$ of intersections of hyperplanes in $\mathcal{A}$ (viewed as affine subspaces of $k^n$), partially-ordered by reverse inclusion. If 
$\bigcap_{H \in \mathcal A} H \ne \emptyset$ then the intersection poset is a geometric lattice, otherwise it's a geometric semilattice.

If $\mathcal{A}$ is a real hyperplane arrangement, then $\R^n\setminus \mathcal{A}$ is disconnected. By the number of regions of $\mathcal A$ we mean the number of  connected components of $\R^n\setminus \mathcal{A}$.  This number can be detected solely from $\mathcal{L}(A)$ as Zaslavsky's formula (\ref{zaseq}) shows.

If $\mathcal{A}$ is a complex hyperplane arrangement, its complement $M_{\mathcal A}:=\C^n\setminus \bigcup_{H \in \mathcal A} H$ is a manifold whose Betti numbers $\beta_i$ can be detected solely from $\mathcal{L}(\mathcal A)$.  Indeed, this follows from the  formula of Orlik and Solomon \cite[Theorem 5.2]{Orlik_Solomon},
\begin{equation}\label{OSeq} \sum_{i =0}^n \beta_i(M_{\mathcal A}) t^i=(-t)^r \chi_{\mathcal L(\mathcal A)}(-t^{-1}),
\end{equation}
where $r$ is the length of $\mathcal L(\mathcal A)$.

\vspace{.1in} \noindent {\bf  The Bond Lattice of a Graph.} 
Let $G$ be a graph on vertex set $[n]$. The {\em bond lattice of $G$} is the subposet $\Pi_G$ of the partition lattice $\Pi_n$ consisting of partitions $\pi = B_1|\cdots|B_k$ such that $G|_{B_i}$ is connected for all $i$.  
Note that $\Pi_n$ is the bond lattice of the complete graph $K_n$.  Another example is given below.

\begin{center}
\includegraphics{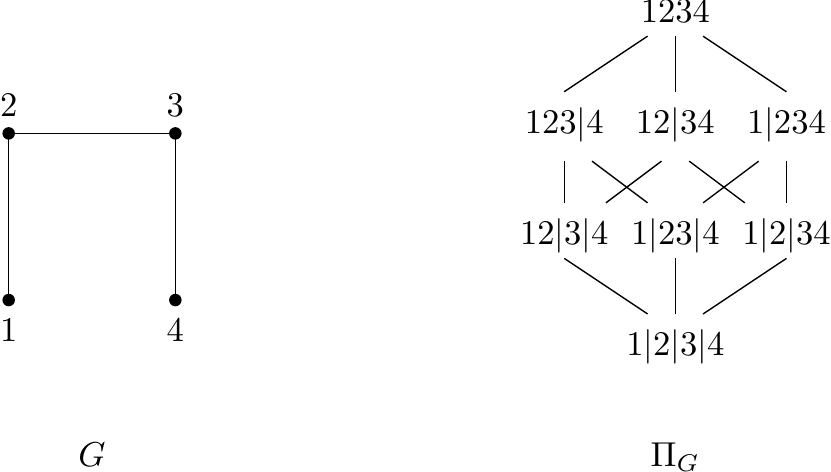}
\end{center}

Broken circuits provide a useful means of computing the M\"obius function of the bond lattice of a graph (or more generally, of geometric lattices).  
Let $G = ([n],E)$ be a finite graph. Fix a total ordering of $E$ and let $S$ be a subset of $E$. Then $S $  is called a {\em broken circuit} if it consists of a cycle in $G$ with its least edge (with respect to this ordering) removed. If $S$ does not contain a broken circuit, we say that $S$ is a {\em non-broken circuit} or {\em NBC} set.

Given any  $S\subseteq E$, let $\pi_S$ be the partition of $[n]$ whose blocks are the vertex sets of the connected components of the graph $([n],S)$. 
The following formula is due to Whitney
 \cite[Section 7]{Logical_Math}  
 for graphs and Rota 
 \cite[Pg. 359]{Found_Comb_1} 
 for general geometric lattices.
For $\pi \in \Pi_G$,
\begin{equation} \label{RWeq} (-1)^{\text{rk}(\pi)}\mu(\hat{0},\pi) = \#\{\mbox{NBC sets } S \mbox { of } G :  \pi_S = \pi\}.\end{equation}

Given a rooted tree whose vertex set is a subset of $\Z^+$, we say the tree is {\em increasing} if each nonroot vertex  is larger than its parent. A rooted forest on a subset of $\Z^+$ is said to be {\em increasing} if it consists of increasing rooted trees.  Note  that if $G$ is $K_n$ then by ordering the edges lexicographically with the smallest element as the first component, the NBC sets of $G$ are exactly the edge sets of the increasing forests on $[n]$.

\vspace{.1in} \noindent
{\bf   Genocchi and median Genocchi Numbers.} The Genocchi numbers and median Genocchi numbers are classical sequences of numbers that have been extensively studied in combinatorics.  There are many ways to define them.  Here we define them in terms of Dumont permutations; see  \cite[p. 44]{Derangements_Genocchi}. 
A {\em Dumont permutation} is a permutation $\sigma \in \mathfrak{S}_{2n}$ such that $2i>\sigma(2i) $ and $2i-1 \le \sigma(2i-1)$ for all $i=1,\dots,n$. A {\em Dumont derangement} is a Dumont permutation without fixed points, i.e.,  $2i > \sigma(2i) $ and $2i-1 < \sigma(2i-1)$ for all $i=1,\dots,n$.

\begin{example}
When $n=2$, the Dumont permutations on $[4]$ (in cycle form) are $$(1,2)(3,4) \qquad (1,3,4,2) \qquad (1,4,2) (3).$$ 
When $n=3$, the Dumont derangements on $[6]$ are:
\begin{center}$\begin{array}{cccc}
(1,3,5,6,4,2) & (1,3,4,2)(5,6) & (1,2)(3,4)(5,6) & (1,2)(3,5,6,4)\\
(1,4,3,5,6,2) & (1,5,6,3,4,2) & (1,5,6,2)(3,4) & (1,4,2)(3,5,6) .
\end{array}$\end{center}
\end{example}

For  $n \ge 1$, the (signless) {\em Genocchi number} $g_n$ is defined to be the number of Dumont permutations on $[2n-2]$,  and for $n \ge 0$, the {\em median Genocchi number} $h_n$ is defined to be the number of Dumont derangements on $[2n+2]$.  We have,

\begin{center}\begin{tabular}{|c|c|c|c|c|c|c|c|}
\hline{\color{red} $n$}  & {\color{red}0} & {\color{red}1}& {\color{red}2}& {\color{red}3} & {\color{red}4}& {\color{red}5} & {\color{red}6}
\\ \hline\hline 
$g_n$ & &$1$ & $1$ & $3$ & $17$ & $155$ & $2073$ 
\\   \hline 
$h_n$ &  $1$ & $2$& $8$ & $56$ & $608$ & $9440$ & $198272$
\\ \hline 
\end{tabular} \end{center}

%===============================================================================

\section{The (type A) homogenized Linial arrangement}
 In this section we give a characterization of the intersection lattice $\mathcal L(\mathcal H_{2n-1})$ as an induced subposet of $\Pi_{2n}$ and compute its M\"obius function.  

\subsection{The intersection lattice is a bond lattice} \label{bondsec} We begin by showing that $\mathcal L(\mathcal H_{2n-1})$ is isomorphic to  the bond lattice of a nice bipartite graph.
Let $\Gamma_{2n}$ be the bipartite graph\footnote{The graph $\Gamma_{2n}$ belongs to a class of graphs called Ferrers graphs, which were introduced by Ehrenborg and van Willegenburg \cite{Enumerative_Ferrers_Graphs} and further studied in \cite{Boolean_Complex_Ferrers}, \cite{Selig_Smith_Steingrimsson}. We have been able to extend some of our results to more general Ferrers graphs and to skew Ferrers graphs. }
 on vertex set $\{1,3,\dots,2n-1\} \sqcup \{2,4,\dots,2n\}$ with an edge between $2i - 1$ and $2j$ iff $i\leq j$.
The graph $\Gamma_6$ is given below.

\begin{figure}[H]
\begin{center}
 \includegraphics{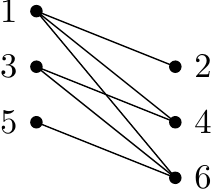}
\end{center}
\end{figure}

\vspace{-.4in} \begin{theorem} \label{bondth} The intersection lattice  $\mathcal L(\mathcal H_{2n-1})$ is isomorphic to the bond lattice $\Pi_{\Gamma_{2n}}$, which is  the induced subposet of $\Pi_{2n}$ consisting of the partitions $X = B_1|\cdots|B_k$ in which $\min(B_i)$ is odd and $\max(B_i)$ is even for all nonsingleton $B_i$. \end{theorem}
In \cite{Full_Paper}, we prove Theorem~\ref{bondth} by constructing an invertible $\Z$-linear map that sends $\mathcal H_{2n-1}$ to an  arrangement whose intersection poset is $\Pi_{\Gamma_{2n}}$.

We  use the Rota-Whitney formula (\ref{RWeq})  to compute the M\"obius function of $\Pi_{\Gamma_{2n}}$.  Our NBC sets have a nice description which we give now.  We  say that a rooted forest on node set $A\subset \Z^+$ is {\em increasing-decreasing (ID)} if the trees are rooted at their largest node and for each internal node $a$,
\begin{itemize}
\item if $a$ is odd then $a$ is less than all its descendants and all its children are even,
\vspace{-.1in} \item if $a$ is even then $a$ is greater than all its descendants and all its children are odd. \end{itemize}

\begin{center}
\includegraphics[height=1.2in]{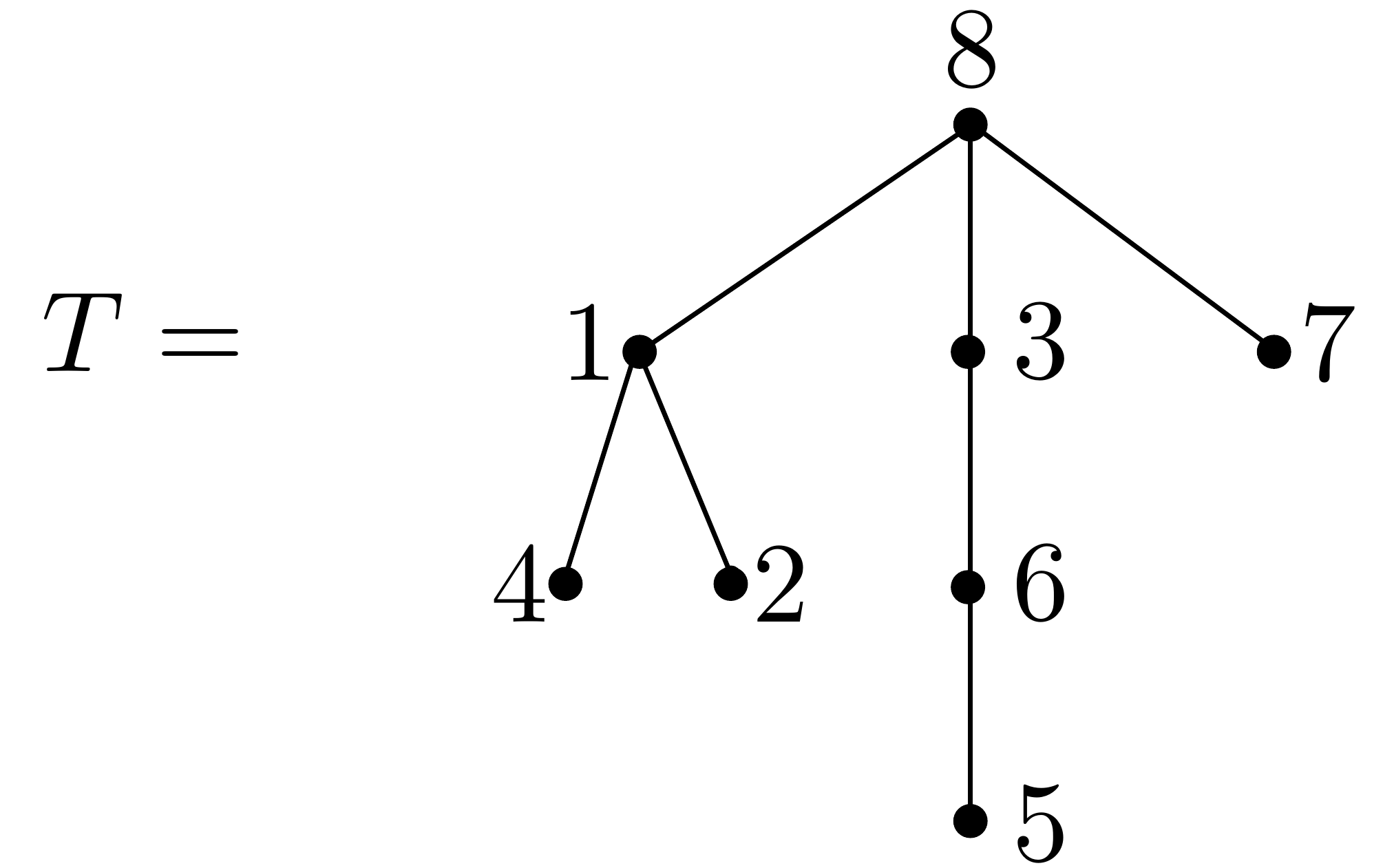}
\end{center}

\begin{theorem}\label{doubth} For all $\pi \in \Pi_{\Gamma_{2n}}$, we have that $(-1)^{|\pi|} \mu(\hat 0,\pi)$ equals the number of ID forests on $[2n]$ whose trees have nodes sets equal to the blocks of $\pi$. Consequently
$$\chi_{\mathcal H_{2n-1}}(t) = \sum_{k=0}^{2n} (-1)^{k} f_{2n,k} t^{k-1} .$$
 where $f_{2n,k}$ is the number of ID forests on $[2n]$ with exactly $k$ trees.
\end{theorem}
To prove this, we show that under an appropriate ordering of the edges of $\Gamma_n$, the NBC sets of $\Gamma_{2n}$ are  the edge sets of the ID forests on $[2n]$.\footnote{This connection between NBC sets and ID forests is also implicit in the proof of Theorem 7.3 of \cite{Selig_Smith_Steingrimsson}.}

\subsection{Dumont-like permutations}

Our next step is to introduce a class of permutations similar to the Dumont permutations and then give a bijection between these permutations in $\s_{2n}$ and  the ID forests on $[2n]$.

Let $A$ be a finite subset of $\Z^+$.  We say $\sigma \in \s_A$ is a {\em D-permutation} on $A$ if 
$i \le \sigma(i) $ whenever $i$ is odd and $i \ge \sigma(i)$ whenever $i$ is even.
We denote by $\mathcal{D}_A$ the set of D-permutations on $A$ and by $\mathcal{DC}_A$ the set of D-cycles on $A$. If $A = [n]$,  we write $\mathcal{D}_n$ and $\mathcal{DC}_{n}$.

Note that all Dumont permutations are D-permutations, but not conversely. Indeed, the only difference between the two classes of 
permutations on $\s_{2n}$ is that fixed points can be even or odd in a D-permutation, while only odd fixed points 
are allowed in a Dumont permutation.  It follows immediately from the definitions that
$$\mathcal{DC}_{2n} \subseteq \{\mbox{Dumont derangements in } \s_{2n}\} \subseteq  \{\mbox{Dumont 
permutations in } \s_{2n}\} \subseteq \D_{2n}.$$  Recall that the two sets in the middle of this chain are 
enumerated by median Genocchi number $h_{n-1}$ and Genocchi number $g_{n+1}$, respectively.  It turns out that the sets on the ends of the chain are also enumerated by Genocchi and median Genocchi numbers. Indeed, there is a bijection between Dumont permutations and D-cycles that  yields, $$|\mathcal{DC}_{2n}| = g_{n},$$   for all $n \ge 1$. We use the theory of surjective pistols discussed in \cite{Une_Famille} and \cite{Derangements_Genocchi} to prove the following result.

\begin{theorem} 
\label{genth}
For all $n \ge 0$,
$|\mathcal{D}_{2n}| = h_{n}$.

\end{theorem}

Next we define a bijection $\psi_A$ from the set $\mathcal T_A$  of  ID trees on $A$ to $\mathcal{DC}_{A}$.
 For $T \in \mathcal T_A$,  order the children of each even node in increasing order and the children of each odd node in decreasing order.  This turns $T$ into a rooted planar tree, which can be traversed in postorder.   Let $\alpha:=\alpha_1,\cdots,
\alpha_{|A|}$ be the word obtained by traversing $T$ in postorder, that is, $\alpha_i$ is the $i$th node of $T$ in postorder.  Now let $\psi_A(T)$ be the permutation whose cycle form is $(\alpha)$.  For the ID tree $T$ given  in Section~\ref{bondsec}, we have $\psi_{[8]}(T)= (4,2,1,5,6,3,7,8)$.

 \begin{lemma} \label{bijlem} For all $A \subseteq [2n]$, the map $\psi_A:\mathcal T_A \to \mathcal{DC}_{A}$ is a well-defined bijection.  Consequently $|\mathcal T_A| = |\mathcal{DC}_{A}|$.
 \end{lemma}

The {\em cycle support} of $\sigma \in \s_n$ is the partition $\cyc(\sigma) \in \Pi_n$ whose blocks are comprised of the elements of the cycles of $\sigma$.  For example,
$\cyc((1,7,2,4)(5)(6,8,9,3)) = 1247|5|3689$.
As a consequence of Theorem~\ref{doubth} and Lemma~\ref{bijlem}, we have the following result.
\begin{theorem} \label{mobth}
For  $\pi \in  \Pi_{\Gamma_{2n}}$, where $n \ge 1$,
$$(-1)^{|\pi|}\mu_{\Pi_{\Gamma_{2n}}}(\hat{0},\pi) = |\{\sigma \in \mathcal{D}_{2n} \st \cyc(\sigma) = \pi\} |.$$ Consequently, \begin{equation} \label{chardumeq} \chi_{\mathcal L(\mathcal H_{2n-1})}(t) = \sum_{k=1}^{2n} s_D(2n,k) t^{k-1},\end{equation}
where $(-1)^{k}s_D(2n,k)$ is equal to  the number of D-permutations on $[2n]$ with exactly $k$ cycles.
\end{theorem}

Next we invoke  Theorem~\ref{genth}.  By setting $t = -1$ in (\ref{chardumeq})  we  recover Hetyei's result  (\ref{heteq}) that $\chi_{\mathcal L(\mathcal H_{2n-1})}(-1)= -h_{n}$ , and by setting $t=0$ we obtain the following new result on the Genocchi numbers.\footnote{For another class of 
 induced subposets of the partition lattice 
 whose  M\"obius invariant can be expressed in terms of the Genocchi numbers, see \cite{sund}.} 

\begin{theorem} \label{muth} For all $n \ge 1$,
$$\mu_{\mathcal L(\mathcal H_{2n-1})}(\hat{0},\hat{1}) = -|\mathcal{DC}_{2n}| = -g_{n}.$$
\end{theorem}

In the full version of the paper \cite{Full_Paper}, we use (\ref{chardumeq}) and the theory of surjective pistols in \cite{Une_Famille} to derive the following result. 

\begin{theorem} \begin{equation} \label{geneq} \sum_{n\geq 1} \chi_{ \mathcal L(\mathcal H_{2n-1})}(t) \, x^n =  \sum_{n\geq 1}\frac{  (t-1)_{n} (t-1)_{n-1} \,x^n}{\prod_{k=1}^n(1-k(t-k)x)},\end{equation}
where $(a)_{n}$ denotes the falling factorial  $a(a-1)\cdots (a-n+1)$.
\end{theorem}
Equation (\ref{geneq}) reduces to a formula of Barsky and Dumont \cite[Lemma 2]{Barsky_Dumont} for the Genocchi numbers when $t$ is set equal to $0$ and to the formula  of Barsky and Dumont  for the median Genocchi numbers given in (\ref{BDeq}) when $t$ is set equal to $-1$.
 
We are also able to obtain  the following characterization of the median Genocchi numbers by evaluating $\chi_{\Pi_{\Gamma_{2n}}}(t)$ in another way.

\begin{theorem} For all $n \ge 1$, $h_{n}$ is equal to the number of permutations $\sigma$ on $[2n]$ whose descents  $\sigma(i) >\sigma(i+1)$ occur only when $\sigma(i)$ is even and $\sigma(i+1)$ is odd.
\end{theorem}

We now give yet another way in which the median Genocchi numbers arise. 
\begin{theorem}\footnote{This result was also independently conjectured by Hetyei (personal communication).}
 \label{affineth} For all $n \ge 3$, 
$$ \chi_{ \mathcal L(\mathcal H_{2n-1})}(t) = (t-1)^3\chi_{P_{n}}(t),$$
where $P_{n} $  is the  intersection semilattice of a certain affine hyperplane arrangement in $\R^{2n-4}$.  Moreover, $|\chi_{P_n}(1)| = h_{n-3}$; hence the number of bounded regions of this affine arrangement is $h_{n-3}$.
\end{theorem}

In the full version of the paper, we compute  $\chi_{P_n}(1)$ by applying the theory of shellability to the NBC complex of the geometric semilattice $P_n$.  The consequence follows from  Zaslavsky's result on the number of bounded regions of an affine arrangement \cite{Facing_Up_Arrangements}.

%=================================================================================================

\section{The homogenized Linial-Dowling arrangement}
In this section, we extend the results of the previous section to the Dowling arrangements, which generalize the complexified  types A and B braid arrangements.

The {\em Dowling lattice} $Q_n(\Z_m)$ consists of labeled partitions $B_0 | B_1 | \dots | B_k$ of $\{0\} \cup [n]$ such that 
\begin{itemize}
\item $0 \in B_0$ ($B_0$ is  called the {\em zero block}),
\item the elements of $B_i$, $i \ge 1$, are labeled with elements of $\{0,1,\dots,m-1\}$ and $\min(B_i)$ is labeled with $0$.
\end{itemize}
The  cover relation is given by merging blocks as follows. Let $B_0|B_1|\cdots|B_k \in Q_n(\Z_m)$. 
\begin{itemize}
\item If $B_0$ and $B_i$ merge, erase all of the labels from $B_i$ and merge the blocks as in $\Pi_n$ to obtain a new zero block $B_0'$.
\item Suppose $i,j \neq 0$, and $\min(B_i) < \min(B_j)$. There are $m$ ways to merge $B_i$ and $B_j$.  For each $\ell \in \{0,\dots,m-1\}$, when  $B_i$ and $B_j$ merge, the labels of the elements of $B_i$ remain unchanged, while $\ell$ is added mod $m$ to the labels of the elements of $B_j$.\end{itemize}

\begin{example}
Suppose $m=3$.  Then $05|1^03^1|2^04^2$ is covered by
$$0135 | 2^04^2 \qquad 0245 | 1^03^1 \qquad 05|1^02^03^14^2 \qquad 05|1^02^13^14^0 \qquad 05|1^02^23^14^1 .$$
\end{example}

It is not hard to see that for all $m \ge 1$, the Dowling lattice $Q_n(\Z_m)$ is isomorphic to the intersection lattice of the Dowling arrangement defined in (\ref{doweq}).  
See 
\cite{Dowling_Geometric_Lattices} and 
\cite{Michelle_Dowling} 
for further information on Dowling lattices.

Now we introduce a Dowling analog of the homogenized Linial arrangement.
Let $\omega = e^{2\pi i/m}$ be a primitive $m$th root of unity. The {\em homogenized Linial-Dowling arrangement} is the complex hyperplane arrangement
$$\mathcal H^m_{2n-1} = \{x_i - \omega^{\ell}x_j = y_i \st 1\leq i < j \leq n, \,\, 0\leq \ell < m\} \cup \{x_i = y_i \st 1\le i \le n\}\subseteq \C^{2n}.$$
Note that when $m=2$, the arrangement $\mathcal{H}^m_{2n-1}$ is a complexified version of the type B homogenized Linial arrangement $\mathcal H^B_{2n-1}$ defined in the introduction.  When $m=1$,  the arrangement $\mathcal{H}^m_{2n-1}$ is the complexified version of the arrangement obtained by  intersecting $\mathcal{H}_{2n-1}$ with the coordinate hyperplane $x_{n+1}=0$.  The resulting arrangement on the coordinate hyperplane has the same intersection lattice as $\mathcal{H}_{2n-1}$.

The proof of the following result is similar to that of the type A version.
 \begin{theorem}
\label{Dowbondth}
For all $n,m \ge 1$, the intersection lattice $ \mathcal L(\mathcal H^m_{2n-1})$  is isomorphic to the induced subposet $\mathcal L_{2n-1}^m$ of $Q_{2n-1}(\Z_m)$ consisting of all labeled partitions such that
\begin{itemize}
\item  for nonsingleton $B_0$, $\min(B_0\setminus\{0\})$ is odd,
\item for all nonsingleton $B_i$, with $i>0$, $\min(B_i)$ is odd and $\max(B_i)$ is even.
\end{itemize}
\end{theorem}

To compute the M\"obius function of the geometric lattice $ \mathcal L(\mathcal H^m_{2n-1})$, we apply the Rota-Whitney formula (\ref{RWeq}) to $\mathcal L_{2n-1}^m$ and then we construct a bijection from the NBC sets of  $\mathcal L_{2n-1}^m$  to the class of $m$-labeled D-permutations, which we now define.  

An $m$-labeled D-permutation $\sigma$ on $[2n]$  is a D-permutation  whose entries are decorated with elements of $\{0,1,\dots,m-1\}$ such that
\begin{itemize}
\item cycle minima are labeled $0$,
\item if $(a_1,a_2, \dots, a_r=2n)$ is the cycle of $\sigma$ containing $2n$ then all  right-to left minima of the word $a_1a_2\dots a_r$ are labeled $0$. \end{itemize}
For example, let $n=5$ and $\sigma = (3,7,8,5,9,10) (1,4,2) (6)$.  Since the right-to-left minima of the first cycle are $10,9,5,3$, they must all be labeled $0$.  Since $1$ and $6$ are the minima of their respective cycles, they must also be labeled $0$.  Hence  $\sigma$ with the labeling  $(3^0,7^*,8^*,5^0,9^0,10^0) (1^0,4^*,2^*) (6^0)$,  where $*$ denotes any label in $\{0,1,2\}$,  is a $3$-labeled D-permutation.

We write $\mathcal{D}_{2n}^m$ for the set of $m$-labeled D-permutations on $[2n]$ 
and $\mathcal{DC}_{2n}^m$ for the set of $m$-labeled D-cycles on $[2n]$. 
 The {\em cycle support} of $\sigma \in \mathcal{D}_{2n}^m$ is the $m$-labeled partition $\cyc(\sigma) = B_0|\cdots|B_k\in Q_n(\Z_m)$ obtained from $\sigma$ as follows:
\begin{itemize}
\item The set of entries of the cycle of $\sigma$  that contains $2n$ gives rise to the zero block $B_0$, with 2n replaced by $0$ and all labels removed.
\item Each cycle of $\sigma$ that doesn't contain $2n$  gives rise to a labeled block $B$ for which the labels of the entries of $B$ are the same as the labels of the entries of the cycle.
\end{itemize}
For example,  if $\sigma =(1^{0}3^{1}4^{1}2^{2})(5^{0})(6^0)(7^{0}8^{0})$ then  
$\cyc(\sigma) = 07|1^02^23^14^1|5^0|6^0$.

\begin{theorem}
 Let $n,m \ge 1$.  For all $\pi  \in \mathcal L_{2n-1}^m$, 
$$(-1)^{|\pi|}\mu_{ \mathcal L_n^m}(\hat{0},\pi) = |\{\sigma \in \mathcal{D}_{2n}^m \st \cyc(\sigma) = \pi\}|.$$
Consequently,
\begin{equation} \label{mchareq}\chi_{ \mathcal L(\mathcal H^m_{2n-1})}(t) = \sum_{k=1}^{2n} s_{D,m}(2n,k) t^{k-1},\end{equation}
where $(-1)^{k}s_{D,m}(2n,k)$ is equal to   the number of $\sigma \in \mathcal D_{2n}^m$ with exactly $k$ cycles, 
and $$\chi_{ \mathcal L(\mathcal H^m_{2n-1})}(-1) = - |\mathcal D_{2n}^m|
 \qquad  \mbox{ and }\qquad  \mu_{\mathcal L(\mathcal H^m_{2n-1})}(\hat{0},\hat{1})= -|\mathcal {DC}_{2n}^m|.$$ 
\end{theorem}

By   the Orlik-Solomon formula (\ref{OSeq}), equation (\ref{mchareq}) gives  a combinatorial formula for  the the Betti numbers of the complement of $\mathcal H_{2n-1}^m$ in $\C^{2n}$.  

In the full version of the paper \cite{Full_Paper}, we use (\ref{mchareq}) and the theory of surjective pistols in \cite{Une_Famille} to derive the following $m$-analog of  (\ref{geneq}).
Note that this reduces to 
(\ref{BBDeq}) when we set $m=2$ and $t=-1$.

\begin{theorem} For all $m \ge 1$, 
\begin{equation} \label{mgeneq} \sum_{n\geq 1}\chi_{ \mathcal L(\mathcal H^m_{2n-1})}(t) \, x^n =  \sum_{n\geq 1}\frac{  (t-1)_{n,m} (t-m)_{n-1,m} \,x^n}{\prod_{k=1}^n(1-mk(t-mk)x)}.\end{equation}
where $(a)_{n,m} = a(a-m)(a-2m)\cdots(a-(n-1)m)$.
\end{theorem}

 There is a well-studied polynomial analog of the Genocchi numbers known as the Gandhi polynomials; see \cite[Section 3]{Derangements_Genocchi}. They are defined by   $G_1(x) = x^2$ and $G_{n}(x) = x^2(G_{n-1}(x+1) - G_{n-1}(x))$ for $n \ge 2$. We obtain the following m-analog of Theorem~\ref{muth}.
 \begin{theorem}  For all $n \ge 1$, $$\mu_{\mathcal L(\mathcal H^m_{2n-1})}(\hat{0},\hat{1}) =- m^{2n-1} G_n(m^{-1}).$$
 \end{theorem}

We also obtain an $m$-analog of Theorem~\ref{affineth}, in which the intersection semilattice $P^m_{n} $ of a certain affine arrangement in $\C^{2n-4}$ satisfies 
 $$\chi_{ \mathcal L(\mathcal H^m_{2n-1}) }(t)= (t-m)(t-1)^2 \chi_{P_n^m}(t).$$

\acknowledgements{The authors thank Jos\'{e} Samper for a valuable suggestion pertaining to bipartite graphs.  They also thank the referees for their comments.  MW thanks G\'abor Hetyei for introducing her to his work on this topic when she visited him at the University of North Carolina, Charlotte during her Hurricane Irma evacuation, and  for his hospitality.
   }

\bibliographystyle{amsplain}
\bibliography{HLBiblio}

\providecommand{\bysame}{\leavevmode\hbox to3em{\hrulefill}\thinspace}
\providecommand{\MR}{\relax\ifhmode\unskip\space\fi MR }
% \MRhref is called by the amsart/book/proc definition of \MR.
\providecommand{\MRhref}[2]{%
  \href{http://www.ams.org/mathscinet-getitem?mr=#1}{#2}
}
\providecommand{\href}[2]{#2}
\begin{thebibliography}{10}

\bibitem{Finite_Field_Method}
C.~A. Athanasiadis, \emph{Characteristic polynomials of subspace arrangements
  and finite fields}, Adv. Math. \textbf{122} (1996), no.~2, 193--233.

\bibitem{Barsky_Dumont}
D.~Barsky and D.~Dumont, \emph{Congruences pour les nombres de {G}enocchi de 2e
  esp\`ece}, Study {G}roup on {U}ltrametric {A}nalysis. 7th--8th years:
  1979--1981 ({P}aris, 1979/1981) ({F}rench), Secr\'{e}tariat Math., Paris,
  1981, pp.~Exp. No. 34, 13.

\bibitem{Boolean_Complex_Ferrers}
A.~Claesson, S.~Kitaev, K.~Ragnarsson, and B.~E. Tenner, \emph{Boolean
  complexes for {F}errers graphs}, Australas. J. Combin. \textbf{48} (2010),
  159--173.

\bibitem{Dowling_Geometric_Lattices}
T.~A. Dowling, \emph{A class of geometric lattices based on finite groups}, J.
  Combinatorial Theory Ser. B \textbf{14} (1973), 61--86.

\bibitem{Interpretations_Combinatoires}
D.~Dumont, \emph{Interpr\'{e}tations combinatoires des nombres de {G}enocchi},
  Duke Math. J. \textbf{41} (1974), 305--318.

\bibitem{Derangements_Genocchi}
D.~Dumont and A.~Randrianarivony, \emph{D\'erangements et nombres de
  {G}enocchi}, Discrete Math. \textbf{132} (1994), no.~1-3, 37--49.

\bibitem{Enumerative_Ferrers_Graphs}
R.~Ehrenborg and S.~van Willigenburg, \emph{Enumerative properties of {F}errers
  graphs}, Discrete Comput. Geom. \textbf{32} (2004), no.~4, 481--492.

\bibitem{Michelle_Dowling}
E.~Gottlieb and M.~L. Wachs, \emph{Cohomology of {D}owling lattices and {L}ie
  (super)algebras}, Adv. in Appl. Math. \textbf{24} (2000), no.~4, 301--336.

\bibitem{Alternation_Acyclic}
G.~{Hetyei}, \emph{{Alternation acyclic tournaments}}, ArXiv e-prints (2017).

\bibitem{Full_Paper}
A.~Lazar and M.~L. Wachs, \emph{On the intersection lattice of the homogenized
  linial arrangement}, In preparation.

\bibitem{Orlik_Solomon}
P.~Orlik and L.~Solomon, \emph{Combinatorics and topology of complements of
  hyperplanes}, Invent. Math. \textbf{56} (1980), no.~2, 167--189.

\bibitem{Defcox}
A.~Postnikov and R.~P. Stanley, \emph{Deformations of {C}oxeter hyperplane
  arrangements}, J. Combin. Theory Ser. A \textbf{91} (2000), no.~1-2,
  544--597, In memory of Gian-Carlo Rota.

\bibitem{Une_Famille}
A.~Randrianarivony and J.~Zeng, \emph{Une famille de polyn\^omes qui interpole
  plusieurs suites classiques de nombres}, Adv. in Appl. Math. \textbf{17}
  (1996), no.~1, 1--26.

\bibitem{Found_Comb_1}
G-C. Rota, \emph{On the foundations of combinatorial theory. {I}. {T}heory of
  {M}\"{o}bius functions}, Z. Wahrscheinlichkeitstheorie und Verw. Gebiete
  \textbf{2} (1964), 340--368 (1964).

\bibitem{Selig_Smith_Steingrimsson}
T.~{Selig}, J.~P. {Smith}, and E.~{Steingrimsson}, \emph{{EW-tableaux,
  Le-tableaux, tree-like tableaux and the Abelian sandpile model}}, ArXiv
  e-prints (2017).

\bibitem{EC1}
R.~P. Stanley, \emph{Enumerative combinatorics. {V}olume 1}, second ed.,
  Cambridge Studies in Advanced Mathematics, vol.~49, Cambridge University
  Press, Cambridge, 2012.

\bibitem{sund}
S.~Sundaram, \emph{The homology of partitions with an even number of blocks},
  J. Algebraic Combin. \textbf{4} (1995), no.~1, 69--92. \MR{1314560}

\bibitem{Logical_Math}
H.~Whitney, \emph{A logical expansion in mathematics}, Bull. Amer. Math. Soc.
  \textbf{38} (1932), no.~8, 572--579.

\bibitem{Facing_Up_Arrangements}
T.~Zaslavsky, \emph{Facing up to arrangements: face-count formulas for
  partitions of space by hyperplanes}, Mem. Amer. Math. Soc. \textbf{1} (1975),
  no.~issue 1, 154, vii+102.

\end{thebibliography}

\end{document}